\documentclass[12pt]{amsart}
\usepackage{amsfonts}
\usepackage{amsmath}
\usepackage{amscd}
\usepackage{amssymb}
\usepackage[dvips]{graphicx}
\usepackage[dvips]{epsfig}
\usepackage{color}

\def\?{\char'76}
\def\!{\char'74}
\def\R{\mathbb{R}}

\def\X{\widetilde{X}}
\def\s{\mathbb{S}^1}  %
\def\Im{\text{\rm Im}}
\def\codim{\text{\rm codim}}
\def\L{\mathcal{L}}

\def\F{\mathfrak{F}}

\def\X{\widetilde{X}}
\def\strat#1{\mathcal{S}_{#1}}
\def\stratmin#1{\strat{#1}^{^{min}}}
\def\M{\Sigma^{^{min}}}
\def\Iso#1{\text{\rm Iso}(#1,\strat{#1})}
\def\Isoeq#1#2{\text{\rm Iso}_{#2}(#1,\strat{#1})}
\def\larrow#1{\overset{\L_{#1}}\longrightarrow}

\newcounter{numero}
\newcommand{\Numero}{\setcounter{numero}{1}(\arabic{numero}) }
\newcommand{\numero}{\addtocounter{numero}{1}(\arabic{numero}) }

\newcounter{letra}
\newcommand{\Letra}{\medskip \setcounter{letra}{1}(\alph{letra}) }
\newcommand{\letra}{\medskip \addtocounter{letra}{1}(\alph{letra}) }

\newcounter{romano}
\newcommand{\Romano}{\setcounter{romano}{1}(\roman{romano}) }
\newcommand{\romano}{\addtocounter{romano}{1}(\roman{romano}) }

\newtheorem{Teo}{Theorem}[section] % en la conf. original es [subsection]
\newtheorem{Lema}[Teo]{Lemma}
\newtheorem{Propo}[Teo]{Proposition}
\newtheorem{Corollary}[Teo]{Corollary}
\newtheorem{conj}[Teo]{Conjecture}

\theoremstyle{definition}

\newtheorem{ejem}[Teo]{Example}
\newtheorem{ejems}[Teo]{Examples}
\newtheorem{falsa}[Teo]{}

\theoremstyle{remark}
\newtheorem{obs}[Teo]{$\clubsuit$ Remark}
\newtheorem{nota}[Teo]{\large{$\dag\dag$\tt Beware}}
\newtheorem{comentarios}[Teo]{Comments}

\def\bteo{\begin{Teo}}
\def\eteo{\end{Teo}}
\def\blema{\begin{Lema}}
\def\elema{\end{Lema}}
\def\bprop{\begin{Propo}}
\def\eprop{\end{Propo}}
\def\bcor{\begin{Corollary}}
\def\ecor{\end{Corollary}}
\def\bconj{\begin{conj}}
\def\econj{\end{conj}}
\def\bfalsa{\begin{falsa}}
\def\efalsa{\end{falsa}}
\def\bejem{\begin{ejem}}
\def\eejem{\end{ejem}}
\def\bejems{\begin{ejems}}
\def\eejems{\end{ejems}}
\def\bobs{\begin{obs}}
\def\eobs{\end{obs}}
\def\bnota{\begin{nota}}
\def\enota{\end{nota}}
\def\bcom{\begin{comentarios}}
\def\ecom{\end{comentarios}}
\def\bdem{\begin{proof}}
\def\edem{\end{proof}}

\begin{document}

\pagestyle{myheadings} \markboth{F. Dalmagro}{Equivariant Unfoldings
for G-Stratified Pseudomanifolds}

\title{Equivariant Unfoldings of G-Stratified Pseudomanifolds}
\author{Dalmagro F.}
\address{Escuela de Matematicas, Facultad de Ciencias, UCV, Caracas
Venezuela}
\email{dalmagro@euler.ciens.ucv.ve}
\dedicatory{To Rodolfo Ricabarra, in memoriam.}
\date{July 2003}
\keywords{Intersection Cohomology, Stratified Pseudomanifolds}
\subjclass{35S35; 55N33}

\maketitle

\begin{abstract}
    For any abelian compact Lie group $G$, we introduce a family
    of $G$-stratified pseudomanifolds, whose main feature is the preservation
    of the orbit spaces in the category of stratified pseudomanifolds.
    Which generalize a previous definition given in \cite{popper}.
    We also find a sufficient condition for the existence of equivariant unfoldings,
    so we have Intersection Cohomology with differential
    forms, as defined in \cite{S}.
    Moreover, if $G$  act on a manifold $M$, we find
    a equivariant  unfolding of $M$ which induce a
    canonical unfolding on the $k$-orbits space for every closed subgroup
    $K$ of $G$.   
\end{abstract}

\section*{Introduction}

Let $X$ be a Thom-Mather stratified space with depth $d(X)=n$.
The De Rhan Intersection Cohomology
 of $X$ with diferential forms was defined in
\cite{brylinsky} by means of an auxiliar construction called {\it unfolding}, which is
a continuous map $\L:\X\rightarrow X$
where $\widetilde{X}$ is a smooth manifold obtained trough a finite composition
\[
    \L:\X = X_n\larrow{n} X_{n-1}\rightarrow\cdots
    \rightarrow X_{1}\larrow{1} X_{0}=X
\]
of topological operations $X_i\larrow{i} X_{i-1}$, called {\it elementary unfoldings}.
This iterative construction is possible because the stratification of
$X$ is controlled by the existence of a
family of conical fiber bundles over the singular strata.
Later in \cite{S} we find a more abstract definition of unfoldings,
which impose some conditions of transversality over the singular strata.
For instance, if the depth of $X$ is  1 then the first elementary unfolding of $X$
is an unfolding in the new sense. \newline

Now let $G$ be a compact Lie group. We introduce the definition of a
$G$-stratified pseudomanifold  in the category of stratified
pseudomanifolds.
Our definition is related to a previous one given in \cite{popper}.
A $G$-stratified pseudomanifold is a stratified pseudomanifold in the usual sense
together with a continuous action preserving the strata, and
whose local model near each singular strata is given by a conical slice.
We also give a definition of equivariant unfoldings, which is a suitable
adaptation of the usual definition of 
unfolding to the family of $G$-stratified pseudomanifolds.
We give a sufficient condition for the existence of equivariant unfoldings, which
is related to the choice of a good family of tubular neighborhoods and a sequence
of equivariant elementary unfoldings of $X$. Each
elementary unfolding induces an elementary unfolding on each orbit space
together with a factorization diagram.\newline

The content of this paper is as follows:

In \S1 we introduce the category of stratified pseudomanofolds.

In \S2 we define stratified $G$-pseudomanifols and study their corresponding
$K$-orbit spaces, for a closed subgroup $K$ of $G$.  

In \S3 we study equivariant tubular neighborhoods, which are equivariant
versions of the usual ones.

In \S4 we introduce the concept of an equivariant explosion of a
$G$-pseudomanifold, and present elementary explosions as given in
\cite{BHS}.

In \S5 we define $G$-pseudomanifolds $X$ stratified with a transversal
Thom-Mather structure whose main feature is that the sequence of elementary
explosions determine an equivariant explosion of the space $X$, which
project naturally onto an equivariant explosion of $X/K$, for a closed
subgroup $K$ of $G$.

\section{Stratified Pseudomanifolds}
In this section we review the usual definitions of stratified spaces,
stratified morphisms and stratified pseudomanifolds. For a more detailed
introduction see \cite{borel2}, \cite{pflaum}.
\bfalsa\label{def espacios estratificados}
    {\bf Stratified spaces}
    Let $X$ be a Hausdorff, locally compact and 2nd countable space.
    A {\bf stratification} of $X$ is a locally finite partition $\strat{X}$ satisfying:
    \item[\Romano] Each element $S\in \strat{X}$ is a connected
    manifold with the induced topology, which a {\bf stratum} of $X$.
    \item[\romano] If $S'\cap\overline{S}\neq\emptyset$ then $S'\subset\overline{S}$
    for any two strata $S,S'\in\strat{X}$. In this case we write $S'\leq S$ and we
    say that $S$ {\bf incides} on $S'$.\\
    We say that $(X,\strat{X})$ is a {\bf stratified space} whenever
    $\strat{X}$ is a stratification of $X$.
\efalsa

With the above conditions, the incidence relationship is a partial
order on $\strat{X}$. More over, since $\strat{X}$ is
locally finite, any strictly ordered chain
\[
        S_0<S_1<\dots<S_m
\]
in $\strat{X}$ is finite. The {\bf depth} of $X$ is by definition
the supremum (possibly infinite) of the integers $m$ such that there is a
strictly ordered chain as above. We write this as $d(X)$.\newline

The maximal (resp. minimal) strata in $X$ are open (resp. closed) in $X$.
A {\bf singular} stratum is a non-maximal stratum in $X$. The  union of the singular strata
is the {\bf singular part} of  $X$, denoted by $\Sigma\subset X$, which is closed in $X$.
Its complement $X-\Sigma$ is open and dense in $X$. The family of minimal strata
will often be denoted by $\stratmin{X}$, while the union of minimal strata
will be denoted by $\M$, which we call the {\bf minimal part} of $X$.

\bfalsa\label{ejems espacios estratificados}{\bf Examples}
    Here there are some examples of stratified spaces.
    \item[\Numero] For any manifold $M$ the trivial stratification of $M$
    is the family
    \[
        \strat{M}=\{C: C\text{  is a connected component of }M\}
    \]

    \item[\numero] For any connected manifold $M$, the space $M\times X$ is
    a estratified space, with the estratification
    \[
	\strat{M\times X} = \{ M\times S: S\in\strat{X}\}  
    \]
    Notice that $d(M\times X)=d(X)$.
    \item[\numero] The {\bf cone} of a compact stratified space $L$ is 
    the quotient space
    \[
        c(L)=L\times [0,\infty)/L\times\{0\}
    \]
    We write $[p,r]$ for the equivalence class of $(p,r)\in L\times[0,\infty)$.
    The symbol $*$ will be used for the equivalence class of
    $L\times\{0\}$, this is the {\bf vertex} of the cone.
    The family
    \[
        \strat{c(L)}=\{*\}\cup\{S\times(0,\infty):S\in\strat{L}\}
    \]
    is the canonical stratification of $c(L)$. Notice that $d(c(L))=d(L)+1$.
\efalsa
\bfalsa{\bf Stratified subspaces and morphisms} Let $(X,\strat{X})$ be a stratified
    space. For each subset $Z\subset X$ the {\bf induced partition} is the family
    \[
        \strat{Z/Y}=\{C:C \text{ is a connected component of } Z\cap S,
	S\in\strat{X} \}
    \]
    We will say that $Z$ is a {\bf stratified subspace} of $X$, whenever the induced
    partition on $Z$ is a stratification of $Z$.\newline

    Now let $(Y,\strat{Y})$ be another stratified space.
    A {\bf morphism} (resp. {\bf isomorphism}) is a
    continuous map $f:X\rightarrow Y$ (resp. homeomorphism) which smoothly
    (resp. diffeomorphicaly) sends strata into strata.
    In particular, $f$ is a {\bf embedding} if $f(X)$ is a stratified
    subspace of $Y$ and
    $f:X\rightarrow  f(X)$ is an isomorphism.
\efalsa

Henceforth, we will write $\Iso{X}$ for the group of isomorphisms of a stratified space $X$.
The following statement will be used later, we leave the proof to the reader.

\blema\label{lema espacios estratificados} Let $(X,\strat{X})$ be a stratified
    space, and
    %\item[\Numero] For each $f\in \Iso{X}$ and each $S\in \strat{X}$, if
    %$S$ and $f(S)$ are comparable strata then $f(S)=S$.
    $\F\subset\strat{X}$ a subfamily of equidimensional strata.
    The connected components of  $M=\underset{S\in\F}{U}S$
    are the strata in $\F$.
\end{Lema}

Stratified pseudomanifolds were used by Goresky and MacPherson in order to introduce
the Intersection Homology and extend the
Poincar\'e duality to the family of stratified spaces.
For a brief introduction the reader can see \cite{borel2}.
\bfalsa\label{def psve}{\bf Stratified pseudomanifolds}
    The definition of a stratified pseudomanifold is made by induction on the
    depth of the space. More precisely:
    \item[\Numero] A stratified pseudomanifold of depth 0 is a manifold
    with the trivial stratification.
    \item[\numero] An arbitrary stratified space $(X,\strat{X})$
    is a {\bf stratified pseudomanifold} if, for any singular stratum $S\in \strat{X}$,
    there is a compact stratified pseudomanifold $L_S$ depending on $S$ (called
    a {\bf link} of $S$)  such that each point $x\in S$ has a coordinate neighborhood
    $U\subset S$ and an embedding onto an open subset of $X$.
    \[
        \varphi:U\times c(L_S)\rightarrow X
    \]
    such that $x\in \Im(\varphi)$. The pair $(U,\varphi)$ is a {\bf chart} of $x$
    {\bf modelled} on $L_S$.
\efalsa
\bfalsa\label{ejems pseudovariedades}
    {\bf Examples} Here are some examples of stratified pseudomanifolds.
    \item[\Numero] If $X$ is a stratified pseudomanifold, then any open subset $A\subset X$
    is also a stratified pseudomanifold. Also the product $M\times X$
    (with the canonical stratification) is a stratified pseudomanifold,
    for any manifold $M$.

    \item[\numero] If $L$ is a compact stratified pseudomanifold, then $c(L)$ is a stratified
    pseudomanifold.
\efalsa

%----------------------------G-pseudovariedades estratificadas------------------------------
\section{$G$-Stratified Pseudomanifolds}
From now on, we fix an abelian compact Lie group $G$.
We will study the family of actions of $G$ which preserve the strata.
Our definition is strongly related to the previous one given in \cite{popper}.
Also some easy proofs in this section can be seen in \cite{gysin1}.\newline

Given a
stratified space $(X,\mathcal{S}_{X})$
and a effective action $\Phi:G\times X\rightarrow X$;
we write $\Phi(g,x)=gx$ for any $g\in G$, $x\in X$. We denoted
$X/K$ by the $K$-orbit space for every $K$ closed subgroup of $G$,
and by $\pi:X\rightarrow X/K$ the orbit map.
The group of $G$-equivariant isomorphisms of $X$ will be denoted by $\Isoeq{X}{G}$.

\bfalsa\label{def G espacio estratificado}
    {\bf $G$-stratified spaces}
    We say that $X$ is {\bf $G$-stratified} whenever:
    \item[\Numero] For each stratum $S\in\strat{X}$ the points of
    $S$ all have the same isotropy group, denoted by $G_S$.
    \item[\numero] Each $g\in G$ induces an isomorphism
    $\Phi_g:X\rightarrow X \in\Isoeq{X}{G}$.\\
    The orbit space $X/K$ inherits a {\bf canonical stratification} given by the
    family
    \[
        \strat{X/K}=\{\pi(S):S\in\strat{X}\}
    \]
    Notice also that $d(X)=d(X/K)$.
\efalsa
\bfalsa\label{ejemstra}
    {\bf Examples} Here are some examples of $G$-stratified spaces:
    \item[\Numero] Each $G$-manifold $M$ has a natural structure of $G$-stratified
    space, when $M$ is endowed with the stratification given by orbit types.
    \item[\numero] If $X$ is a $G$-stratified space, then
    $M\times X$ is a $G$-stratified space with the action $g(m,x)=(m,gx)$; for any manifold $M$.
    \item[\numero] If $L$ is a compact $G$-stratified space then $c(L)$, with the action
    $g[x,r]=[gx,r]$, is a $G$-stratified space.
\efalsa

\blema\label{lema accion estratificada cociente}
    Let $G$ be an abelian compact Lie group, $K\subset G$ a closed subgroup.
    Then for any $G$-stratified space $X$ the orbit space $X/K$ is a
    $G/K$-stratified space.
\elema
\bdem
    %First notice that $X$ is also a $K$-stratified space and $X/K$
    %inherits a canonical stratification
    %induced by the orbit map $\pi:X\rightarrow X/K$.
    Write $\overline{g}\in G/K$ for the equivalence class of $g\in G$.
    Consider the quotient action
    \[
        \overline{\Phi}:G/K\times X/K\rightarrow X/K
        \hskip1cm
	\overline{g}\cdot\pi(x)=\pi(gx)
    \]
    This action is well defined because $G$ is abelian. So:\\
    $\bullet$ {\tt The isotropy groups are constant over the strata of $X/K$:}
    This is straightforward, since for each stratum $S\in\strat{X}$ we have
    \[
        (G/K)_{\pi(S)}=KG_S/K
    \]
    Hence $\pi(S)$ has constant isotropy.\\
    $\bullet$ {\tt Each $\overline{g}$ induces an isomorphism
    $\overline{\Phi}_g\in G/K \in \Iso{X/K}$:}
    For each $g\in G$ we have a $K$-equivariant isomorphism $\Phi_g\in\Iso{X}$.
    Passing to the quotients we obtain an isomorphism
    $\overline{\Phi}_g\in G/K \in \Iso{X/K}$.
    The differentiability of this map on $\pi(S)$ is immediate from the
    following commutative diagram
    \[
        \begin{CD}
        S @>{\Phi_g}>> gS \\
        @V{\pi} VV @VV{\pi} V \\ 
        \pi(S) @>{\overline{\Phi}_g}>> \pi(gS)\\
        \end{CD}
    \]
\edem

Now we introduce the definition of a $G$-stratified pseudomanifold.
\bfalsa\label{def slice conico y Gpsve} {\bf $G$-stratified pseudomanifolds}
    A $G$-stratified pseudomanifold is a stratified pseudomanifold in the
    usual sense, endowed with a structure of $G$-stratified space (i.e.
    $G$ acts  by isomorphisms) and whose local model is described through
    conical slices. Conical slices were introduced in \cite{popper}
    in order to state a sufficient condition on any continuous action of 
    a compact Lie group (abelian or not) on stratified pseudomanifold 
    so that the corresponding orbit space would
    remain in the same class of spaces.\newline

    Let $(X,\strat{X})$ be a $G$-stratified space.
    Take a singular stratum $S\in\strat{X}$ a point $x\in S$.
    A {\bf conical slice} of $x$ in $X$ is a slice $S_x$ in the usual sense of \cite{bredon},
    with a conical part transverse to the stratum $S$. In other words:
    \item[\Numero] $S_x$ is an invariant $G_S$-space containing $x$.
    \item[\numero] For any $g\in G$, if $gS_x\cap S_x\neq \emptyset$
    then $g\in G_S$.
    \item[\numero] $GS_x$ is open in $X$. And
    \item[\numero] There is a $G_S$-equivalence $\beta: \R^i\times c(L)\rightarrow S_x$
    where $i\geq 0$ and $L$ is a compact $G_S$-stratified space.
    Here the action of $G_S$ on $\R^i$ is trivial
    (notice that $\beta$ induces on $S_x$ a structure of
    $G_S$-stratified space).\newline

    The definition of a {\bf $G$-stratified pseudomanifold} is made by
    induction on the depth of the space.
    A $G$-stratified pseudomanifold with depth 0 is a manifold with a smooth
    free action of $G$. In general, we will say that $X$ is a $G$-stratified 
    pseudomanifold
    if, for each singular stratum $S\in\strat{X}$,
    there is a compact $G_S$-stratified pseudomanifold
    $L_S$ such that each point $x\in S$ has a conical slice
    \[
        \beta:\R^i\times c(L_S)\rightarrow S_x
    \]
    and the usual map on the twisted product
    \[
        \alpha:G\times_{G_S}S_{x}\rightarrow X
        \hskip1cm
	\alpha([g,y])=gy
    \]
    is an  equivariant (stratified) embedding on an open subset of $X$.
    We say that the triple $(S_x,\beta,L_S)$ is a {\bf distinguished slice} of $x$.
\efalsa
\bfalsa\label{ejems G-pseudovariedades estratificadas}{\bf Examples}
    Here there are some examples of $G$-stratified pseudomanifolds.
    \item[\Numero] Take a smooth effective action 
    $\Phi:G\times M\rightarrow M$ with fixed points 
    on a manifold $M$ endowed with the
    stratification by orbit types. By the
    Equivariant Slice Theorem, $M$ is a $G$-stratified pseudomanifold.
    \item[\numero] If $X$ is a $G$-stratified pseudomanifold then
    $M\times X$ is a $G$-stratified pseudomanifold with the obvious action.
    \item[\numero] If $L$ is a compact $G$-stratified pseudomanifold, then
    $c(L)$ is a $G$-stratified pseudomanifold with the obvious action.
    \item[\numero] Any invariant open subspace of a $G$-stratified pseudomanifold
    is itself a $G$-stratified pseudomanifold.
    %\item[\numero] Any Thom-Mather space $X$ endowed with an action of $G$
    %preserving the tubular neighborhoods is a $G$-stratified pseudomanifold.
\efalsa

\bobs\label{Prop1.2.5}
    {\it Each $G$-stratified pseudomanifold is a
    stratified pseudomanifold in the previous sense.}
    
    To see this, proceed by induction on the depth.
    Take a $G$-stratified pseudomanifold $X$.
    For $d(X)=0$ the statement is trivial.
    Assume the inductive hipothesis and suppose that $d(X)>0$. Take a singular
    stratum $S\in\strat{X}$, a point $x\in S$ and a distinguished slice
    $(S_x,\beta,L_S)$ of $x$.
    The isotropy subgroup $G_S$ acts on $G$ by the restriction of the group operation.
    We fix a slice $S_e$ of the identity element $e\in G$ with respect to this action.
    Since $G_S S_e$ is open in $G$, the composition
    {\small\begin{eqnarray*}
        (S_{e}\times \R^{i})\times c(L_S) \rightarrow
        S_{e}\times (\R^{i}\times cL_{S})\rightarrow
	S_{e}\times S_{x}\rightarrow\\
	\rightarrow S_{e}\times (G_S\times_{G_S}S_{x})\rightarrow
	(G_S S_{e})\times_{G_S}S_{x}\rightarrow X
    \end{eqnarray*}}
    is an embedding. Notice that $L_S$ is a stratified pseudomanifold by induction.
    Since $S_e\times \R^{i}\simeq S_e G_S(S\cap
    S_x)$ is open in $S$. We have obtained a chart of $x$ modelled on $L_S$.
\eobs
\bobs\label{obs X es K-psve}
    {\it If $X$ is a $G$-stratified pseudomanifold and $K$ is any
    closed subgroup of $G$, then $X$ is also a $K$-stratified
    pseudomanifold.}  

    It is straightforward that $X$ is a $K$-stratified
    space. For any singular stratum $S$ and any $x\in S$, in order to choose
    a distinguished slice in $x$ we proceed as follows: Take a distinguished
    slice $\beta:\R^i\times c(L_S)\rightarrow S_x$ in $x$ with respect to the
    action of $G$. Take also a slice $V_e$ of the identity element $e\in G$
    with respect to the action of $G_SK$ in $G$. Then
    $\imath\times\beta:(V\times \R^i)\times c(L_S)\rightarrow VS_x$  is a distinguished
    slice of $x$ with respect to the action of $K$.
\eobs

Now we study the factorization of a $G$-stratified pseudomanifold
when considered as a $K$-stratified pseudomanifold for any closed
subgroup $K\subset G$.

\bprop\label{prop X es G/K-psve}
    Let $G$ be a compact, abelian Lie group;
    $K\subset G$ a closed subgroup.
    If  $X$ is a $G$-stratified pseudomanifold
    then $X/K$ is a $G/K$-stratified pseudomanifold.
\eprop
\bdem
    As before, write $\pi:X\rightarrow X/K$ for the orbit map induced by
    the action of $K$ on $X$.
    Proceed by induction on $l=d(X)$. For $l=0$
    it is straightforward, since $d(X/K)=d(X)=0$.
    Assume the inductive hipothesis and suppose that $d(X)>0$.
    By \S\ref{lema accion estratificada cociente},
    $X/K$ is a $G/K$-stratified space, so we must verify the existence
    of conical slices.\newline

    Take a singular stratum $S\in\strat{X}$, fix a point $x\in S$
    and a distinguished slice $(S_x,\beta,L_{S})$ of a $x$.
    The $K$-equivariant isomorphism $\beta:\R^i\times c(L_S)\rightarrow S_x$
    induces an isomorphism on the orbit spaces
    \[
        \overline{\beta}:\R^i\times c(L_S/G_S\cap K)\rightarrow \pi(S_x)
        \hskip1cm
        \overline{\beta}(b,[\overline{l},r])=\pi(\beta(b,[l,r]))
    \]
    Now we will show that the triple $(\pi(S_x),\overline{\beta},L_S/G_S\cap K)$
    is a distinguished slice of $\pi(x)\in X/K$. We do it in three steps.
    \item[$\bullet$] {\tt $\pi(S_x)$ is a slice of $\pi(x)$:} This is straightforward,
    since $(G/K)_{\pi(x)}=KG_S/K$, the quotient $\pi(S_x)$ is a
    $(G/K)_{\pi(x)}$-space with the quotient action and
    the orbit map $\pi$ is an open map.
    \item[$\bullet$] {\tt $\overline{\beta}$ is a $KG_S/K$-equivalence:}
    This is immediate, since $\beta$ is an $H$-equivalence.
    Notice that, by induction on the depth,
    $L_{S}/G_S\cap K$ is a $KG_S/K$-stratified pseodumanifold.
    \item[$\bullet$] {\tt The induced map $\overline{\alpha}:
    (G/K)\times_{(G_S/G_S\cap K)}\pi(S_x)\rightarrow X/K$ is an embedding:}
    This $\overline{\alpha}$ is given by the rule
    $\overline{\alpha}([\overline{g},\pi(z)])=\overline{g}.\pi(z)$, and
    is a homeomorphism. We consider the following commutative diagram
    \[
    \begin{CD}
    G\times_{G_S}S_x @>\alpha >> X \\
    @V{\overline{\pi}} VV @VV\pi V \\
    (G/K)\times_{\tiny(G_S/G_S\cap K)}\pi(S_x) @>\overline{\alpha}>> X/K\\
    \end{CD}
   \]
    Since the vertical arrows are submersions, and $\alpha$ is an embedding, we obtain
    that $\overline{\alpha}$ is an embedding.
\edem

%-------------------------------Entornos tubulares-----------------------------------

\section{Tubular neighborhoods}

Henceforth we fix a compact, abelian Lie group $G$, and a $G$-stratified pseudomanifold
$X$. In this section we will study the family of equivariant tubular neighborhoods,
which are  equivariant version of the usual ones. Given a singular stratum $S$ in $X$,
a tubular neighborhood is just a locally trivial fiber  bundel over a $S$,
whose fiber is $c(L_S)$, the cone of the link of $S$; and whose structure group is
$\Isoeq{L_S}{G_S}$. \newline

We will clarify these ideas inmediately, considering a previous step
in our way: the definition of a $G$-stratified fiber bundle. The reader will find in
\cite{steenrod} a detailed introduction to the fiber bundles, while \cite{thom} provides
the usual definition of a tubular neighborhood in the stratified context (see also
\cite{bredon} for the smooth case).

\bfalsa\label{def G-fibrado estratificado}
    {\bf $G$-stratified fiber bundles}
    Let $\xi=(E,p,B,F)$ be a locally trivial fiber bundle with (maximal)
    trivializing atlas $\mathcal{A}$. We will say that $\xi$ is a
    {\bf $G$-stratified} whenever:
    \item[\Numero] The total space $E$ is a $G$-stratified space.
    \item[\numero] The base space $B$ is a manifold, endowed with a smooth action
    $\Psi:G\times B\rightarrow B$ and with constant isotropy $H\subset G$
    at all its  points.
    \item[\numero] The fiber $F$ is a $H$-stratified space.
    \item[\numero] The projection $p:E\rightarrow B$ is $G$-equivariant.
    \item[\numero] The group $G$ acts by isomorphisms. In other words,
    each chart
    \[
        \varphi:U\times F\rightarrow p^{-1}(U) \in\mathcal{A}
    \] is $H-$equivariant;
    and for any two charts $(U,\varphi),(U',\varphi^{'})\in\mathcal{A}$
    such that $U'\cap g^{-1}U\neq\emptyset$ for some $g\in G$, there is
    a map
    \[
        g_{\varphi,\varphi^{'}}:U'\cap g^{-1}U\rightarrow \Isoeq{F}{H}
    \]
    such that
    \[
             \varphi^{-1}g\varphi'(b,z)=(gb,g_{\varphi,
             \varphi^{'}}(b)z)
    \]
\efalsa

\blema\label{lema G-fibrados estratificados}
    Let $\xi=(E,p,B,F)$ be a  $G$-stratified fiber bundle,
    $H$ the isotropy of $B$.
    If $F$ is an $H$-stratified pseudomanifold, then
    $E$ is a $G$-stratified pseudomanifold.
\elema

\bdem
    Fix a singular stratum $S$ in $E$  and a point $x\in S$. We must prove
    the existence of a link $L_S$ depending only on $S$ and,
    a distinguished slice $(S_x,\beta,L_S)$ in
    $x$. For this purpose, let's take  a trivializing chart
    \[
        \varphi:U\times F\rightarrow p^{-1}(U)\in\mathcal{A}
    \]
    such that $x\in p^{-1}(U)$. Take $z=p(x)$ and a $G$-slice $V_{z}$ in
    $B$. Since $V_z$ is contractible, we assume that
    $V_z\cong\R^k$ and $V_z$  is contained in $U$. \newline

    Write $\varphi^{-1}(x)=(z,y)\in V_z\times F$ and take $S'$ the stratum in $F$
    containing $y$. Since $F$ is an $H$-stratified pseudomanifold, we can choose
    a distinguished slice $S_y$ in $y$; say
    \[
        \beta_0:S_y\rightarrow\R^i\times c(L_{S'})
    \]
    Consider the following composition
    \[
        \varphi(V_{z}\times S_{y})
	\overset{\varphi^{-1}}\rightarrow
	V_z\times S_y
	\overset{\imath\times\beta_0}\rightarrow
	V_z\times\R^i\times c(L_{S'})
	\cong R^{i+k}\times c(L_{S'})
    \]
    We will show that
    \[
        (S_x,\beta,L_S)=(\varphi(V_{z}\times S_{y}),(\imath\times\beta_0)\circ\varphi^{-1},
	L_{S'})
    \]
    is a distinguished slice in $x$. We proceed in three steps.\newline

    $\bullet$ {\tt $L_S$ only depends on $S$:} If $(U',\psi)\in\mathcal{A}$
    is another trivializing chart covering $x$,
    $\psi^{-1}(x)=(z,y')\in V_z\times F$ and $\beta'_0:
    S_{y'}\rightarrow \R^i\times c(L_{S''})$ is a distinguished slice in
    $y'$; then the composition  $\beta'\beta^{-1}$ induces an $H$-isomorphism
    $L_S\overset{\cong}\rightarrow L_{S''}$.\newline

    $\bullet$ {\tt $S_x$ is a conical slice:}
    We verify the conditions (1) to (4) of \S\ref{def slice conico y Gpsve}.
    \item[\Numero] Since $V_z$ is a slice of $z\in B$, we have
    $gp(x)=p(gx)=p(x)\in V_z$ for any $g\in G_S$. So $G_S=H\cap G_S=H_S$,
    but $\varphi$ is $H$-equivariant, hence $G_S=H_S=H_{S'}$.
    Again, since $\varphi$ is $H$-equivariant and $S_y$ is $H_{S'}=G_S$
    invariant, we obtain that $S_x$ is $G_S$-invariant.
    \item[\numero] Take $g\in G$, $x'\in S_x$ such that $gx'\in S_x$.
    Then $gp(x')=p(gx')\in V_z$, so $g\in H$ and $gp(x')=p(x')$.
    Since $\varphi$ is $H$-equivariant, if
    $x'=\varphi(p(x'),y)$ then $g.x'=\varphi(p(x'),gy)$, and
    $gy\in S_y$; hence $g\in H_{S'}=G_S$.
    \item[\numero] Take a slice $S_e$ of the identity element $e\in G$ with respect
    to the action of $H$. Since $S_e$ is contractible, we can assume that
    $S_eV_z\subset U$. Notice that $S_eH$
    is open in  $G$. Since $GS_x=\underset{g\in
    G}{\cup}g(S_eH)S_x$,
    we only have to show that $(S_eH)S_x$ is open in $X$. But
    $\varphi$ is $H$-equivariant and the action of $H$ on
    $V_z$ is trivial, so we get the following equality
    \[
            (S_eH)S_x=S_e\left(H\varphi(V_z\times S_y)\right)=
            S_e\varphi(V_z\times HS_y)
    \]
    Since $HS_y$ is open in $F$ we deduce that $S_e\varphi(V_zHS_y)$
    is open in $S_e\varphi(V_z\times F)$.
    Finally we show that $S_e\varphi(V_z\times
    F)=S_ep^{-1}(V_z)$ is open in $X$:
    Since $p$ is equivariant and $S_eV_z$ is open in $U$
    the set $S_ep^{-1}(V_z)=p^{-1}(S_eV_z)=p^{-1}(S_eHV_z)$ is open in
    $p^{-1}(U)$ (and so in $X$).
    \item[\numero] It is straightforward that the map $\beta$ is a $G_S$-equivalence.\newline

    $\bullet$ {\tt $S_x$ is a distinguished slice:}
    We will show that usual
    the map
    \[
        \alpha:G\times_{G_S}S_x\rightarrow X
    \]
    is a (stratified) embedding.
    \item[\Letra] {\tt $\alpha$ preserves the strata:}
    Take a stratum $S^0$ in $S_x$. We will prove that $G'S^0$
    is an open subset in some stratum of $X$, for any connected
    component $G'\subset G$. It is enough to prove this for the connected
    component $G_0$ of the identity element $e\in G$.
    Let $H_0$ be the connected component of the identity element
    $e\in H$. The set $S_eH_0$ is a connected open subset in
    $S_eH$, so is also connected and open in $G_0$. Since
    $G_0S^0$ is connected, we need to prove that $S_eH_0S^0$
    is open in some stratum of $X$. But $S_eHS_x$
    is contained in $p^{-1}(S_eV_z)$ and $\varphi$
    is a stratified embedding, and so we only have to show that
    $\varphi^{-1}(S_eH_0S^0)$ is open in some stratum of $(S_eV_z)\times F$.
    Consider the map
    \[
         \begin{array}{c}
                \Delta :S_eH\times V_z\times S_y \rightarrow
                (S_eV_z)\times F \\
                ^{(gh,b,l) \mapsto
	       (ghb,(gh)_{\varphi\varphi}(b)(z))=(gb,g_{\varphi\varphi}(b)(hz))}
          \end{array}
    \]
    Let $S^1$ be the stratum of $S_y$ such that $S^0=\varphi(V_z\times
            S^1)$. By hypothesis $S_y$ is a distinguished slice of $y$ in $F$, 
	    and there is a stratum $S^2$ in  $F$ such that $H_0S^1$ is open in $S^2$.
            Notice that
	    \[
	         \varphi^{-1}(S_eH_0S^0)=\Delta(S_eH_0\times
                 V_z\times S^1)=
                 \Delta(S_e\times V_z\times H_0S^1)
	    \] 
	    Also, since $\varphi$ is $H$-equivariant,
	    we have 
	    \[
	        p(\varphi^{-1}(S_eH_0S^0))=S_eV_z
	    \]
	    Hence the projection $pr_2:U\times F\rightarrow F$
	    sends $\varphi^{-1}(S_eH_0S^0)$ on some open subset of $S^2$.
	    Notice that $\varphi^{-1}(S_eH_0S^0)$ is connected,
	    so
	    \[
	        pr_2(\varphi^{-1}(S_eH_0S^0))=
                \underset{(g,b)\in S_e\times V_z}{\bigcup}
                g_{\varphi\varphi}(b)(H_0S^1)
            \]
            is a connected subset of $F$. Each  $g_{\varphi\varphi}(b)$
	    is an $H$-equivariant stratified isomorphism;
            hence $g_{\varphi\varphi}(b)(H_0S^1)$ is open is
	    some stratum of $F$ with the same dimension of $S^2$.
	    Since
            $e_{\varphi\varphi}(b)(H_0S^1)=H_0S^1\subset
            S^2$, by \S\ref{lema espacios estratificados}-(2)
            the set $\underset{(g,b)\in S_e\times
            V_z}{\bigcup} g_{\varphi\varphi}(b)(H_0S^1)$ is contained in $S^2$.

    \item[\letra] {\tt $\alpha$ is smooth on each stratum:}
            Since $G\times_{G_x}S_x$ has the quotient stratification induced
	    on  $G\times S_x$ by the action of $H$, the
            stratification of $S_x$ is induced by  $X$ and the action of
	    $G$ is smooth on each stratum of $G\times X$.
            We conclude that the restriction of $\alpha$ to each stratum
	    is smooth.
\edem

\bfalsa\label{def entornos tubulares equivariantes}
    {\bf Equivariant tubular neighborhoods}
    An equivariant tubular neighborhood
    is a conical locally trivial fiber bundle.
    For a detailed introduction the reader can see
    \cite{pflaum}, \cite{thom}.
    In \cite{BHS}, the tubular neighborhoods are used
    in order to show  the existence of an unfolding for
    any manifold endowed with a Thom-Mather structure.
    We will provide an equivariant version of this fact
    for any $G$-stratified pseudomanifold.\newline

    Let $X$ be a $G$-stratified pseudomanifold with $d(X)>0$.
    Let's take a singular stratum $S$ in $X$. An
    {\bf equivariant tubular neighborhood} of $S$ is a $G$-stratified fiber bundle
    $(T_S,\tau_{S},S,c(L_S))$ with (maximal) trivializing atlas $\mathcal{A}$,
    verifying
    \item[\Numero] $T_{S}$ is an open invariant neighborhood of $S$ and
    the inclusion $S\rightarrow T_S$ is a section of
    $\tau_S:T_S\rightarrow S$.
    \item[\numero] $G$ preserves the conical radium:
    For any two charts $(U,\varphi),(U',\varphi^{'})\in\mathcal{A}$
    such that $U'\cap g^{-1}U\neq\emptyset$ for some $g\in G$, there is
    a map
    \[
        g_{\varphi,\varphi^{'}}:U'\cap g^{-1}U\rightarrow \text{Iso}_{G_S}(L_S,\strat{L_S})
    \]
    such that
    \[
             \varphi^{-1}g\varphi'\left(b,[l,r]\right)=\left(gb,[g_{\varphi,
             \varphi^{'}}(b)l,r]\right)
    \]
    This allows us to define a (global) {\bf radium} on $T_S$, as the map
    $\rho_{S}:T_S\rightarrow [0,\infty)$
    satisfying
    \[
        \rho_S(\varphi(z,[l,r]))=r \ \forall (z,[l,r])\in U\times c(L_S);(U,\varphi)\in\mathcal{A}
    \]
    We also define the {\bf radial action} $\delta_S:\R^{+}\times T_S\rightarrow T_S$ as follows
    \[
	\delta_S(r,x)=\varphi(z,[l,rt])
	\ \forall(z,[l,t])\in U\times c(L_S);(U,\varphi)\in\mathcal{A}
        \ \ (\mbox{for } \ x=\varphi(z,[l,t])).
    \]
    We will write $rx$ instead of $\delta_S(r,x)$ in the future.
    These functions satisfy
    \begin{itemize}
       \item[\Letra] $\rho_S(rx)=r\rho_S(x)$ and $\rho_{S}(gx)=\rho_{S}(x)$
       for any $r\in \R^{+}$, $x\in T_S$, $g\in G$.
       \item[\letra] $S\cap \rho^{-1}_{S}(0,\infty)=\emptyset$
       \item[\letra] The radial action commutes with the action of $G$ on $T_{S}$.
    \end{itemize}
\efalsa
\bfalsa\label{def thom-mather}
    {\bf Thom-Mather spaces} (see \cite{verona2}, \cite{verona}):
    A {\bf Thom-Mather $G$-stratified pseudomanifold} 
    is a pair $(X,\mathcal{T})$ where $X$ is a $G$-stratified pseudomanifold
    and $\mathcal{T}=\{T_S:S\in\strat{X}^{^{sing}}\}$ is a family of equivariant
    tubular neighborhoods satisfying
    the following condition:
    \begin{center}
    $T_S\cap T_R\neq\emptyset$ $\Leftrightarrow$ $R\leq S$ or
    $S\leq R$
    \end{center}
    for any two singular strata $R,S$ in $X$. We will usually ommit the family
    $\mathcal{T}$ if there is no possible confusion.
\efalsa

\bfalsa\label{ejems entornos tubulares}
    {\bf Examples} Here are some examples of $G$-stratified
    tubular  neighborhoods.
    \item[\Numero] Following \cite[p.306]{bredon}, for any manifold $M$ endowed
    with a smooth action $\Phi:G\times M\rightarrow M$ there is a riemannian
    metric $\mu$ such that $G$ acts by $\mu$-isometries. By the local properties of
    the exponential map, each singular stratum $S$ of $M$ has a
    smooth $G$-equivariant tubular neighborhood which can be realized as the normal
    fiber bundle $N_\mu(S)$ over $S$ with respect to $\mu$.
    The cocycles
    of this bundle are orthogonal actions. Hence, this tubular neighborhood is
    actually a $G$-stratified tubular neighborhood.
    \item[\numero] If $L$ is a compact $G$-stratified pseudomanifold, the
    map $c(L)\rightarrow\{\star\}$ is a $G$-stratified tubular neighborhood
    of the vertex.
    \item[\numero] If $\xi=(T_S,\tau_S,S,c(L_S))$ is a $G$-stratified
    tubular neighborhood of $S$ in $X$, then $(M\times T_S,\imath_M\times\tau_S,
    M\times S,c(L_S))$ is a $G$-stratified tubular neighborhood of $M\times S$
    in $M\times X$; for any connected manifold $M$.
    \item[\numero] If $f:Y\rightarrow X$ is a $G$-equivariant isomorphism,
    then for any $G$-stratified tubular neighborhood
    $\xi=(T_S,S,\tau_S,c(L_S))$ of a stratum $S$ in $X$; the pull-back
    $f^{*}(\xi)=(f^{-1}(T_S),f^{-1}\tau_S f, f^{-1}(S),c(L_S))$
    is a $G$-stratified tubular neighborhood of $f^{-1}(S)$ in $Y$.
\efalsa

\bprop\label{prop el tubo pasa al cociente}
    Let $X$ be a $G$-stratified pseudomanifold, $K$ a closed subgroup of $G$.
    Write $\pi:X\rightarrow X/K$ for the orbit map induced by the action of $K$.
    Let $\xi=(T_S, \tau_S, S, c(L_S))$
    be an equivariant tubular neighborhood of $S$ in $X$ and write
    \[
        \overline{\tau_S}:\pi(T_S)\rightarrow \pi(S)
    \]
    for the induced quotient map. Then
    $\xi/K=(\pi(T_S), \overline{\tau_S},\pi(S), c(L_S/G_S\cap K))$
    is an equivariant
    tubular neighborhood of $\pi(S)$ in $X/K$.
\eprop
\bdem
   Since $\pi$ is an open map, $\pi(T_S)$ is an open neighborhood
   of $\pi(S)$ in $X/K$.  Also the inclusion $\pi(S)\rightarrow \pi(T_S)$
   is a section of $\overline{\tau_S}:\pi(T_S)\rightarrow\pi(S)$.
   In order to prove that $\xi/K$ is a $G$-stratified tubular neighborhood
   we should first verify that it is a $G$-stratified fiber bundle,
   but the conditions \S\ref{def G-fibrado estratificado}-(1) to (4) are
   straightforward.\newline 
   
    Now we will prove
    \S\ref{def entornos tubulares equivariantes}-(2), which implies
    \S\ref{def G-fibrado estratificado}-(5).
    We will show that the trivializing atlas
    $\mathcal{A}=\{(U,\varphi)\}$ of $\xi$ induces a trivializing atlas
    $\mathcal{A}/K=\{(V,\psi)\}$ of $\xi/K$.
    Write $\pi':L_S\rightarrow L_S/G_S\cap K$ for the orbit map induced by the
    action of $G_S\cap K$ in $L_S$.\newline

    $\bullet$ {\tt Trivializing charts}:
    Take a chart $(U,\varphi)\in\mathcal{A}$
    and a point $x\in U$. Take also a $K$-slice $V$ of $x$ in $S$,
    we assume that $V\subset U$. Since $G_S$ acts trivially on
    $V$ and $KV$ is open in $S$ we deduce that
    \[
        V=V/G_S\cap K=\pi(KV)
    \]
    is open in $\pi(S)$. Since $\varphi$ is $G_S$-equivariant, the function
    \begin{equation}\label{eq carta inducida en tubo cociente}
        \psi:V\times c(L_S/G_S\cap K)\rightarrow\pi(T_S)
	\hskip1cm
        \psi(b,[\pi'(l),r])=\pi(\varphi(b,[l,r]))
    \end{equation}
    is well defined. Moreover, $\psi$ is injective because $G$ acts by isomorphisms
    and $V$ is a $K$-slice in $S$. Notice that $W=KV\cap U$ is open in $U$;
    since $G$ also preserves the radium in $T_S$,
    \[
        \Im(\psi)=\pi(\varphi(W\times c(L_S)))
    \]
    Hence $\Im(\psi)$ is open in $X/K$.
    It is straightforward that $\psi$ sends smoothly strata onto strata, so actually
    $\psi$ is an embedding.\newline

    $\bullet$ {\tt Atlas and cocycles}:
    We consider the family $\mathcal{A}/K=\{V,\psi)\}$ of all the
    pairs $(V,\psi)$ as in (\ref{eq carta inducida en tubo cociente}).
    We will show that $\mathcal{A}/K$ is a trivializing
    atlas of $\xi/K$.
    Take two charts $(V,\psi);(V',\psi')\in\mathcal{A}/K$ respectively induced
    by $(U,\varphi);(U',\varphi')\in\mathcal{A}$.
    Assume that there is some $\overline{g_0}\in G/K$ such that
    $\overline{g_0}^{-1}V\cap V'\neq\phi$;
    so $g^{-1}U\cap U'\neq\phi$ for some $g\in g_0K$. By
    \S\ref{def entornos tubulares equivariantes}-(2), there is a map
    \[
        g_{\varphi\varphi'}:g^{-1}U\cap U'\rightarrow \Isoeq{L_S}{G_S}
    \]
    satisfying
    \[
        g\varphi^{'}(b,[l,r])=\varphi(gb,[g_{\varphi\varphi'}(b)(l),r])
	\hskip1cm
        (b,[l,r])\in (g^{-1}U\cap U')\times c(L_S)
    \]
    Passing to the orbit space $L_S/G_S\cap K$ we obtain the induced map
    \[
        \overline{g_0}_{\psi\psi'}:\overline{g_0}^{-1}V\cap V'\rightarrow
        \Isoeq{L_S/G_S\cap K}{(G_S/G_S\cap K)}
    \]
    satisfying
    {\small \[
        \overline{g_0}\psi'(b,[\pi'(l),r])=\psi(\overline{g_0}b,[
	\overline{g_0}_{\psi\psi'}(b)(\pi'(l)),r]);
	\hskip0.3cm
        (b,[\pi'(l),r])\in (\overline{g_0}^{-1}V\cap V')\times c(L_S/G_S\cap K)
    \]}
    Notice that, by definition, $G/K$ preserves the radium of $\pi(T_S)$.
\edem

%--------------------------Explosiones Equivariantes---------------------------------

\section{Equivariant unfoldings}
An unfolding of a stratified pseudomanifold is an auxiliar
construction which allows us
to define the intersection cohomology from the point of view
of differential forms \cite{BHS},\cite{brylinsky}.
For a detailed introduction to unfoldings, the reader can see
\cite{davis}, \cite{S}.
In this section we introduce equivariant unfoldings, these are
a suitable adaptation of the usual unfoldings to the equivariant category.
We also show that for any
$G$-manifold, considered as a $G$-stratified pseudomanifold,
there is always an equivariant unfolding which induces
a canonical unfolding on the orbit space.

\bfalsa\label{def explo e quivariante}
    {\bf Equivariant unfoldings} Broadly speaking, an un unfolding of a stratified
    pseudomanifold $X$ is a manifold $\X$  and a surjective
    continuous map $\L:\X\rightarrow X$ such that
    $\L^{-1}(X-\Sigma)$ is a union of finitely many disjoint copies of $X-\Sigma$,
    and  which smoothly unfolds the singular part so that the restriction
    $\L:\L^{-1}(S)\rightarrow S$ is a submersion, for any singular stratum $S$.\newline

    As for the usual unfoldings, the definition of an equivariant unfolding is
    made by induction on the depth. Let $X$ be a $G$-stratified pseudomanifold.
    An {\bf equivariant unfolding} of $X$ is a manifold $\X$ together with a smooth
    free action $\widetilde\Phi:G\times \X\rightarrow \X$; a surjective,
    continuous, equivariant map
    \[
        \L:\widetilde{X}\rightarrow X
    \]
    and a family of equivariant unfoldings
    $\{\L_{L_S}:\widetilde{L_S}\rightarrow L_S\}_{S\in\strat{X}^{sing}}$
    where $S$ runs on the singular
    strata of $X$; satisfying:
    \item[\Numero] The restriction $\L:\L^{-1}(X-\Sigma)\rightarrow X-\Sigma$
    is a smooth finite trivial covering.
    \item[\numero] For each singular stratum $S$ and each $x\in S$, there
    is a {\bf liftable modelled chart}, i.e.; a commutative square
    \[
            \begin{CD}
                U\times\widetilde{L_S}\times\R @>{\widetilde\varphi} >> \X \\
                @V \L_c VV @VV\L V \\
                U\times c(L_S) @>\varphi >> X\\
            \end{CD}
    \]
    such that
    \begin{itemize}
	\item[\Letra] $(U,\varphi)$ is a $G_S$-equivariant chart of $x$ modelled
        on $L_S$.
	\item[\letra] $\widetilde{\varphi}$ is a $G_S$-equivariant smooth embedding
	on an open subset of $\X$.
	\item[\letra] The map $\L_c$ is given by the rule
        $\L_c(u,z,t)=(u,[\L_{L_{S}}(z),|t|])$.
    \end{itemize}
    A $G$-stratified pseudomanifold $X$ is said to be {\bf unfoldable} whenever it has an
    equivariant unfolding.
\efalsa

\bfalsa{\bf Examples} Here are some examples of equivariant unfoldings.
    \item[\Numero] For any free smooth action $\Phi:G\times M\rightarrow M$
    the identity $\imath:M\rightarrow M$ is an equivariant unfolding.
    \item[\numero] If $\L:\X\rightarrow X$ is an equivariant unfolding,
    then for any manifold $M$ the product $\imath:M\times \X\rightarrow M\times X$ is also
    an equivariant unfolding.
    \item[\numero] For any equivariant unfolding $\L:\widetilde{L}\rightarrow L$
    over a compact $G$-stratified pseudomanifold $L$, the map
    $\L_c:\widetilde{L}\times\R\rightarrow c(L)$ defined above is also
    an equivariant unfolding.
\efalsa

\bfalsa\label{def explosion primaria}
    {\bf Elementary unfolding of a $G$-stratified pseudomanifold}
    The elementary unfolding of a Thom-Mather space
    is essentially the resolution of singularities given in \cite{davis}
    for the smooth case. This topological operation can be done because the stratification
    is controlled through a family of tubular neighborhoods.
    Under certain conditions, after the iterated composition of finitely
    many elementary unfoldings, one obtains an equivariant unfolding as defined above.
    We follow the exposition of \cite{BHS}.\newline

    Henceforth we fix a Thom-Mather $G$-stratified pseudomanifold $X$,
    a closed (hence minimal) stratum $S$ in $X$ and an equivariant tubular neighborhood
    $(T_S,\tau_S,S,c(L_S))$ of $S$. Define the {\bf unitary sub-bundle} as the
    set $E_{S}=\rho^{-1}_{S}(1)$; this is by construction a $G$-invariant
    stratified subspace of $X$. The restriction $\tau_S:E_S\rightarrow S$ is a
    $G$-stratified fiber bundle with fiber $L_S$.
    Consider the map
    \begin{equation}\label{eq explo tubo}
        \L_{T_{S}}:E_{S}\times\R\rightarrow T_{S}
        \hskip1cm
        \L_{T_{S}}(x,t)=\left\{
        \begin{array}{cc}
            |t|*x & \mbox{ si } t\neq  0 \\
	    \tau_{S}(x) & \mbox{ si } t=0
        \end{array}\right.
    \end{equation}
    Each chart $(U,\varphi)$ in the trivializing atlas provides a
    local description of $\L_{T_{S}}$ through the following commutative square
    \[
	\begin{CD}
	U\times L_{S}\times \R @>\widehat{\varphi}>> E_{S}\times\R \\
       Â    @V 1_{U}\times\L_{C} VV 	@VV\L_{T_{S}}V \\
	U\times cL_{S} @>\varphi >> T_{S}
	\end{CD}
    \]
    where $\widehat{\varphi}(x,l,t)=(\varphi(x,[l,1],t))$ and $\L_{C}(l,t)=[l,|t|]$.
    We also obtain the following properties:
    \item[\Letra] The map $\widehat{\varphi}$ is a $G_{S}$-equivariant embedding.
    \item[\letra] The composition $\tau_{S}\circ \L_{T_{S}}:E_{S}\times\R\rightarrow S$
    is a locally trivial fiber bundle with fiber $L_{S}\times\R$
    and structure group $\Isoeq{L_S}{G_S}$.
    \item[\letra] $d(E_{S}\times\R)=d(E_{S})=d(T_S)-1$. \newline

    Now take a disjoint family of equivariant tubular neighborhoods
    $\{T_S: S\in\stratmin{X}\}$ of the minimal strata.
    The {\bf elementary unfolding} of
    $X$ with respect to the family $\{T_S: S\in\strat{X}^{^{min}}\}$
    is the pair $(\widehat{X},\L)$ constructed as follows:
    First $\widehat{X}$ is the amalgamated sum
    \begin{equation}\label{eq espacio explo primaria}
        \widehat{X}=[\underset{S\in\stratmin{X}}\sqcup E_S\times\R]\
	\underset{\theta}\bigcup\
	[(X-\M)\times\{\pm 1\}]
    \end{equation}
    where $\theta$ is the  map
    \begin{equation}\label{eq pegamento explo primaria}
        \theta:\underset{ S\in\stratmin{X}}\sqcup\ E_S\times\R^*
	\rightarrow
        [X-\M]
	\times\{\pm 1\}
	\hskip1cm
	\theta(x,t)=(|t|*x,|t|^{-1}t)
    \end{equation}
    Second, $\L$ is the continous map given by the rule
    \begin{equation}\label{eq explo primaria}
        \L:\widehat{X}\rightarrow X
	\hskip1cm
        \L(x)=\left\{
	\begin{array}{ll}
	    \L_{T_{S}}(x) & x\in E_S\times\R \\
	    y & x=(y,j)\in(X-\M)\times\{\pm 1\}
        \end{array}
	\right. 
    \end{equation}
\efalsa

Here there are some properties of the elementary unfoldings.

\bprop\label{prop propiedades G-explosiones primarias}
    Let $\L:\widehat{X}\rightarrow X$ be the elementary
    unfolding of a Thom-Mather $G$-stratified pseudomanifold
    $X$. Then
    \item[\Numero] $\widehat{X}$ is a $G$-stratified pseudomanifold, whose
    stratification is the family $\strat{\widehat{X}}$ consisting of
    all the following sets
    \[
        \widehat{R}=
	[ \underset{S\in\stratmin{X}}\sqcup (E_{S}\cap R)\times\R]
        \underset{\theta}{\sqcup} (R\times\{\pm 1\})
    \]
    where $R$ runs over the non closed strata in $X$.
    Moreover, $\widehat{X}$ satisfies the Thom-Mather condition.
    \item[\numero] The map $\L$ is a $G$-equivariant morphism.
    The restriction
    \[
        \L:\L^{-1}(X-\M)\rightarrow
        X-\M
    \]
    is a (trivial) double covering.
    \item[\numero] $d{(\widehat{X})}=d({X})-1$. In particular,
    if $d(X)=1$ then $\L:\widehat{X}\rightarrow X$ is an equivariant
    unfolding.
    \item[\numero] If $X$ is compact, then so is $\widehat{X}$.
    %\item[\numero] $\widehat{M\times X}=M\times \widehat{X}$, for any manifold $M$.
    \item[\numero] For any closed subgroup $K\subset G$, the induced map
    $\overline{\L}:\widehat{X}/K\rightarrow X/K$ is an elementary
    unfolding.
\eprop
\bdem
    \Numero The stratification of $\widehat{X}$ can be seen in \cite{BHS}.
    Since each equivariant tubular neighborhood
    is a $G$-stratified pseudomanifold (because they are invariant
    open subsets of $X$); so are the unitary sub-bundles 
    (see \S\ref{lema G-fibrados estratificados}), and hence $\widehat{X}$
    is a $G$-stratified pseudomanifold. Now we verify the Thom-Mather condition:
    Take a family $\{T_S:S\in\strat{X}\}$ of equivariant tubular neighborhoods
    in $X$. Take also a stratum $\widehat{R}$ in $\widehat{X}$
    induced by a non closed stratum $R$ in $X$. Define
    \[
        T_{\widehat{R}}=\underset{S\in\stratmin{X}}\sqcup\ (E_S\cap T_R)\times\R
	\cup (T_R\times\{\pm 1\})=\L^{-1}(T_R)
    \]
    where $\theta$ is the map given in the equation
    (\ref{eq pegamento explo primaria}) of \S\ref{def explosion primaria}.
    This $T_{\widehat{R}}$ is an equivariant tubular neighborhood of
    $\widehat{R}$ in $\widehat{X}$; we leave the detais to the reader. \\
    \numero and \numero are straightforward, see again \cite{BHS} for more
    details. The last observation of (3) is a consequence of
    def.\S\ref{def explo e quivariante}. \\
    (4) Since $X$ is compact, $\stratmin{X}$ is finite. But
    $\widehat{X}$ is the  quotient of the finite family of compact spaces
    $\underset{S\in\stratmin{X}}{\sqcup}(E_S\times [-1,1])$
    and
    $[X-\underset{S\in\stratmin{X}}{\sqcup}\rho^{-1}_{S}[0,1/2)]\times\{-1,1\}$.
    Then we get the result.\\
   %(5) This is obvious from the definition.\\
   (5) This is a consequence of \S\ref{prop el tubo pasa al cociente}.
\edem

\bobs\label{nota-transversa}
    With tubular neighborhood of \ref{ejems entornos tubulares}-3,
    $\widehat{M\times X}=M\times \widehat{X}$, for any manifold $M$.
    
\eobs

%---------------------------------------------------------------------------

\section{Iteration of elementary unfoldings}\label{def iteracion de explosiones primarias}
    From now on, we fix a Thom-Mather $G$-stratified pseudomanifold $X$.
    We will study the
    composition of finitely many elementary unfoldings, starting at $X$.
    As we have already seen,
    for any elementary unfolding $\L:\widehat{X}\rightarrow X$, the space
    $\widehat{X}$ is again a Thom-Mather $G$-stratified pseudomanifold
    and satisfies $d(\widehat{X})=d(X)-1$. This allows us to ask
    for the behavior of a chain
    \begin{equation}\label{eq cadena explo primaria}
       X_l\larrow{l} X_{l-1}\larrow{l-1} \dots
       \larrow{2} X_1 \larrow{1} X
    \end{equation}
    of elementary unfoldings, where $l=d(X)$. As we shall see, under certain
    conditions on the tubular neighborhoods, this iterative process
    leads us to an equivariant unfolding
    \[
        \L:\X\rightarrow X
    \]
    where $\X=X_l$ and $\L=\L_1\dots\L_{l}$.\newline

    Recall the definition of a {\bf saturated subspace} \cite{BHS}. Let $Y\subset X$
    be a stratified subspace of $X$. We say that $Y$ is {\bf saturated} whenever
    \[
        Y\cap T_S=\tau_S^{-1}(Y\cap S)
	\hskip1cm
	\forall S\in\strat{X}
    \]
    For instance, if $S$ is a singular stratum and $U\subset S$ is open, then
    $Y=\tau_S^{-1}(U)$ is a saturated. Also the unitary sub-bundle
    $Y=E_S$ is saturated.

\bfalsa\label{def morfismos transversos}
    {\bf Transverse morphisms}
    Now we introduce the family of transverse morphisms, whose
    main feature is the preservation of the tubular neighborhoods.
    Let $H\subset G$ be a closed subgroup, $Y$ a Thom-Mather $H$-stratified
    pseudomanifold and $M$ be a connected manifold.
    A morphism
    \[
        \psi: M\times Y\rightarrow X
    \]
    is {\bf transverse} whenever:
    \item[\Numero] $\Im(\psi)$ is a saturated open subspace of $X$.
    \item[\numero] If $\psi(M\times S)\subset R$ then $\psi^{-1}(T_R)=M\times T_S$,
    for any $R\in\strat{X}, S\in\strat{Y}$.\\
    Now let  $\psi:M\times Y\rightarrow X$ be a transverse morphism.
    The  {\bf lifting} of $\psi$ is, by definition, the map
    {\small\[
        \widehat\psi:M\times \widehat{Y}\rightarrow\widehat{X}
        \hskip1cm
        \widehat\psi(m,z,t)=\left\{
        \begin{array}{ll}
           (\psi(m,z),t) & (m,z,t)\in M\times E_S\times\R  \\
	   (\psi(m,z),t) & (m,z,t)\in M\times (Y-\M)\times\{\pm 1\}
        \end{array}
        \right.
    \]}
    This is the unique morphism such that the diagram
    \[
    	\begin{CD}
	M\times \widehat{Y} @>\widehat{\psi}>> \widehat{X} \\
        Â @V {\imath_M\times \L_Y} VV 	@VV \L_X  V \\
	M\times Y @>\psi >> X
	\end{CD}
    \]
    commutes.
\efalsa

\bfalsa\label{ejem morfismo transverso}
    {\bf Examples} For any smooth effective action of $G$ 
    in a manifold, the trivializing
    charts of the tubular neighborhoods are transverse morphisms. In order
    to see this, take a manifold $M$ endowed with a smooth action
    $\Phi:G\times M\rightarrow M$ and an invariant metric $\mu$ in $M$.
    Recall that $M$ has a natural structure of Thom-Mather
    $G$-stratified
    pseudomanifold, where $\strat{M}$ is the stratification induced by
    the orbit types of the action. For any singular stratum $S$ with codimension
    $\codim(S)=q+1>0$, the equivariant tubular neighborhood $T_S=N_\mu(S)$ is the
    normal fiber bundle over $S$ induced by $\mu$
    (see \S\ref{ejems entornos tubulares}). Take also a trivializing chart
    \[
        \varphi: U\times c(\mathbb{S}^q)\rightarrow \tau_S^{-1}(U)
    \]
    We claim that $\varphi$ is transverse. First notice that $\Im(\varphi)$ is a saturated
    open subspace in $M$, so we only have to verify \S\ref{def morfismos transversos}-(2).
    Let $S'$ be a stratum in $c(\mathbb{S}^q)$, $R$ a stratum in $M$. Suppose that
    $\varphi(U\times S')\subset R$. We consider the following cases:\\
    $\bullet$ {\tt $S'=\{\star\}$ is the vertex:} It is straightforward, since $R=S$
    and $T_{S'}=c(\mathbb{S}^q)$.\\
    $\bullet$ {\tt $S'=S''\times\R^{+}$ for some stratum $S''$ in $\mathbb{S}^q$:} Then
    $S<R$. We consider in $T_S$ the following decomposition of the metric:
    \[
        \mu\mid_{T_S}=\mu_H + \mu_V
    \]
    corresponding to the the orthogonal decomposition of the tangent $T(T_S)$ in the
    horizontal and vertical subfiber bundles. Hence
    \[
        \varphi^{-1}(T_R)=\varphi^{*}(N_\mu(R))=N_{\varphi^{*}(\mu)}(U\times S')=
	U\times N_{\mu_V}(S')=U\times T_{S'}
    \]
\efalsa

Now we show two easy properties of the transverse morphisms.
\bprop\label{prop propiedades morfismo transverso}
    Let $K,H$ a closed subgroups of $G$, $L$ a
    Thom-Mather $H$-stratified pseudomanifold,
    $\psi:M\times L\rightarrow X$ a transverse morphism. Then
    \item[\Numero] The lifting $\widehat\psi:M\times\widehat{L}
    \rightarrow \widehat{X}$ is transverse.
    \item[\numero] The induced quotient map
    $\overline{\psi}:M\times (L/H\cap K)\rightarrow X/K$
    is transverse.
\eprop
\bdem
    \Numero is straightforward from the def.\S\ref{def morfismos transversos}.
    \numero is a consequence of \S\ref{prop el tubo pasa al cociente}.
\edem

Finally, we provide a sufficient condition for the existence of an equivariant
unfolding, depending on the transversality of the tubular neighborhoods.
\bteo\label{teo factorizacion explo-primarias}
    Let $X$ be a Thom-Mather $G$-stratified pseudomanifold.
    Suppose that for any singular stratum $S$, each trivializing chart
    \[
        \varphi:U\times c(L_S)\rightarrow T_S
    \]
    is transverse. Then
    \item[\Numero] The composition of the $l$ elementary unfoldings of
    starting at $X$ induces an equivariant unfolding
    $\L:\X\rightarrow X$ where $\X$ is the last (non trivial)
    elementary unfolding and $\L=\L_1\L_2\dots\L_l$ (see eq.
    (\ref{eq cadena explo primaria}) at the begining of this section).
    \item[\numero] For any closed subgroup $K\subset G$, the induced map
    $\overline{\L}:\X/K\rightarrow X/K$ is an unfolding.
\eteo
\bdem
    \Numero Take a family of equivariant tubular neighborhoods in $X$ with transverse
    trivializing charts. Let
    \[
        X_l\larrow{l} X_{l-1}\larrow{l-1} \dots \larrow{2} X_1 \larrow{1} X
    \]
    be the chain of elementary unfoldings induced by this family of tubular
    neighborhoods. We proceed by induction on $l=d(X)$; for $l=1$ the statements
    are trivial. Suppose that $l>1$ and assume the inductive hypothesis,
    so $\L':\X\rightarrow X_1$ is
    an equivariant unfolding, for $\X=X_l$ and $\L'=\L_2\dots\L_l$. Take a closed stratum $S$
    and a transverse trivializing chart
    \[
        \varphi:U\times c(L_S)\rightarrow \tau_S^{-1}(U)\subset T_S
    \]
    Apply the chain of elementary unfoldings and use
    \S\ref{prop propiedades morfismo transverso}; you will get the following
    commutative diagram:
    \[
        \begin{CD}
	    U\times \widetilde{L_S}\times\R @>\widetilde{\psi_1}>> \X \\
            Â @V {\imath_U\times \L_{L_S}\times\imath_{\R}} VV 	@VV \L'  V \\
	    U\times L_S\times\R @> {\psi_1} >> X_1 \\
            Â @V {\imath_U\times \L_S} VV 	@VV \L_1  V \\
	    U\times c(L_S) @>\psi >> X
	\end{CD}
    \]
    We conclude that $\L=\L_1\L':\X\rightarrow X$ is an equivariant unfolding.\\
    \numero This is a consequence of \S\ref{prop propiedades G-explosiones primarias}-(6).
\edem

\bcor[Unfolding of a $G$-manifold]\label{expl}\ \\
    Let $M$ be a manifold, $\Phi:G\times M\rightarrow M$ a smooth effective
    action, possibly with fixed points. Endow $M$ with the stratification
    induced by the orbit types and the usual structure of a Thom-Mather
    $G$-stratified pseudomanifold. Then there is an equivariant unfolding
    $\L:\widetilde{M}\rightarrow M$.
\ecor
\bdem
    Apply the above theorem to the transverse charts obtained in 
    \S\ref{ejem morfismo transverso}.
\edem

\section*{Acknowledgments}
I would like to thank M. Saralegi for some helpful conversations,
and R. Popper for some accurate remarks. While writing this work, I 
received the hospitality of the staff in the Math Department, Universit\'e
D'Artois; and the finnancial support of the ECOS-Nord project
and the CDCH-Universidad Central de Venezuela.

%------------------------------------------------------------------------

%---------------------------------------------------------------------

\end{document}